\documentclass[a4paper,11pt]{amsart}
\usepackage{amsmath,amssymb,colordvi}
\newtheorem{theorem}{Theorem}
\newtheorem{corollary}[theorem]{Corollary}

\newtheorem{example}[theorem]{Example}
\newtheorem{proposition}[theorem]{Proposition}
\newtheorem{definition}[theorem]{Definition}

\newtheorem{remark}[theorem]{Remark}

\newcommand{\blackbox}{\vrule Depth0pt height5pt width5pt}
\newcommand{\bgproof}{\noindent {\bf Proof.} \hspace{2mm}}
\newcommand{\edproof}{\hfill \blackbox \vspace{3mm}}

\begin{document}

\title[Regularity of solutions]
{Regularity of solutions to stochastic Volterra equations with
infinite delay}

\author{Anna Karczewska}
\address{Department of Mathematics, University of Zielona G\'ora,
 ul.\ Szafrana 4a, 65-246 Zielona G\'ora, Poland.}
\email{A.Karczewska@im.uz.zgora.pl}
\thanks{}

\author{Carlos Lizama}
\address{Universidad de Santiago de Chile, Departamento de
Matem\'atica, Facultad de Ciencias, Casilla 307-Correo 2,
Santiago, Chile.}
 \email{clizama@lauca.usach.cl}
\thanks{The second author is supported in part by FONDECYT Grant \#1050084}




\begin{abstract}
In this article we give necessary and sufficient conditions
providing regularity of solutions to stochastic Volterra equations
with infinite delay on a $d$-dimensional torus. The harmonic
analysis techniques and stochastic integration in function spaces
are used. The work applies to both the stochastic heat and wave
equations.
\end{abstract}

\maketitle

\section{Introduction}

We consider the following integro-differential stochastic equation
with infinite delay
\begin{equation}{\label{eq1.1}}
X(t, \theta) = \int_{-\infty}^t b(t-s)[ \Delta X(s, \theta) +
\frac{\partial W_{\Gamma}}{\partial s} (s,\theta)] ds, \quad t
\geq 0, \quad \theta \in T^d,
\end{equation}
where $ b \in L^1 ( \mathbb{R}_+)$, $\Delta$ is the Laplace
operator and $T^d$ is the $d$-dimensional torus. In (\ref{eq1.1}),
$W_{\Gamma}$ is a spatially homogeneous Wiener process with the
space covariance $\Gamma$ and taking values in the space of
tempered distributions $\mathcal{S}'(T^d).$ Such equation arises,
in the deterministic case, in the study of heat flow in materials
of fading memory type (see \cite{ClDaPr88}, \cite{Nu}).

In this paper we address the following question: under what
conditions on the covariance $\Gamma$ the process $X$ takes values
in a Sobolev space $H^\alpha(T^d)$, particularly in $L^2(T^d)$?

This problem has been investigated, for the stochastic heat and
wave equations in case of $ \mathbb{R}^2$, by Dalang and Frangos
\cite{DaFr} (see also Mueller \cite{Mu}). Their results were
recently extended on $\mathbb{R}^d $ and generalized to $T^d$ in
\cite{KaZa2}. In \cite{KaZa3}, regularity of stochastic convolutions
solving Volterra equations on $\mathbb{R}^d$ has been studied.
The results have been obtained by using the resolvent operators
corresponding to Volterra equations and expressed in terms of the
spectral measure and the covariance kernel $\Gamma$ of the Wiener
process $W_\Gamma$. Our approach is different. 
We study a particular case of weak solutions under the basis of an
explicit representation of the solution to (\ref{eq1.1}) (cf.
formula (\ref{eq3.7})).
We remark that the knowledge of the regularity of
solutions is important in the study of nonlinear stochastic
equations (see e.g. \cite{DaFr} and \cite{MSS}).

Observe that equation (\ref{eq1.1}) can be viewed as the {\it
limiting equation } for the stochastic Volterra equation

\begin{equation}{\label{eq1.2}}
X(t, \theta) = \int_{0}^t b(t-s)[ \Delta X(s, \theta) +
\frac{\partial W_{\Gamma}}{\partial s} (s,\theta)] ds, \quad t
\geq 0, \quad \theta \in T^d,
\end{equation}
 see \cite[Chapter III, Section 11.5]{Pr} to obtain details on this assertion. If $b$ is sufficiently regular, we get, by
differentiating (\ref{eq1.2}) with respect to $t$,

\begin{equation}{\label{eq1.3}}
\frac{\partial X}{\partial t}(t, \theta) =
 b(0)[\Delta X(t,\theta)+ \frac{\partial
W_{\Gamma}}{\partial t} (t,\theta)] + \int_0^t b'(t-s)[ \Delta
X(s, \theta) + \frac{\partial W_{\Gamma}}{\partial s} (s,\theta)]
ds,
\end{equation}
where $ t \geq 0$ and $ \theta \in T^d.$\\

Taking in (\ref{eq1.3}), $b(t) \equiv 1$ we obtain
\begin{equation}{\label{eq1.4}}
 \left \{\begin{array}{rcl} \displaystyle
 \frac{\partial X}{\partial t} (t, \theta)
 &=& \Delta X(t,\theta)+ \frac{\partial
W_{\Gamma}}{\partial t} (t,\theta), \quad t >0, \theta \in T^d,
\\ \noalign{\smallskip} X(0,\theta) &= & 0, \quad \theta \in
T^d.
\end{array} \right. \end{equation}
Similarly, taking $b(t) \equiv t$ and differentiating
(\ref{eq1.2}) twice with respect to $t$ we obtain

\begin{equation}{\label{eq1.5}}
 \left \{\begin{array}{rcl} \displaystyle
 \frac{\partial^2 X}{\partial t^2}(t, \theta)
 &=& \Delta X(t,\theta)+ \frac{\partial
W_{\Gamma}}{\partial t} (t,\theta), \quad t >0, \theta \in T^d,
\\ \noalign{\smallskip} X(0,\theta) &= & 0, \quad \theta \in
T^d, \\ \noalign{\smallskip} \frac{\partial X}{\partial t}
(0,\theta) &= & 0, \quad \theta \in T^d.
\end{array} \right. \end{equation}
It has been shown in \cite[Theorem 5.1]{KaZa2} (see also
\cite[Theorem 1]{KaZa1}) that equations
 (\ref{eq1.4}) and (\ref{eq1.5}) on the $d$-dimensional
 torus $ T^d$ have
an $ H^{\alpha + 1}( T^d)$-valued solutions if and only if the
Fourier coefficients $(\gamma_n)$ of the space covariance $\Gamma$
of the process $W_\Gamma$ satisfy
\begin{equation}{\label{eq1.6}}
\sum_{n\in \mathbb{Z}^d} \gamma_n (1 + |n|^2)^{\alpha} < \infty.
\end{equation}

Observe that for both, stochastic heat (\ref{eq1.4}) and wave
(\ref{eq1.5}) equations, the conditions are exactly the same,
despite of the different nature of the equations. On the other
hand, the obtained characterization gives a natural framework in
which nonlinear heat and wave equations can be studied.

In this article, we will prove that condition (\ref{eq1.6}) even
cha\-rac\-te\-ri\-zes  $H^{\alpha + 1}( T^d)$ - valued solutions
for the stochastic Volterra equation (\ref{eq1.1}), provided
certain conditions on the kernel $b$ are satisfied. This is a
strong contrast with the deterministic case, where regularity of
(\ref{eq1.1}) is dependent on the kernel $b$. The conditions that
we impose on $b$ are satisfied by a large class of functions.
Moreover, the important examples $b(t) = e^{-t}$ and $b(t) =
te^{-t}$ are shown to satisfy our assumptions.

We use, instead of resolvent families, a direct approach to the
equation (\ref{eq1.1}) finding an explicit expression for the
solution in terms of the kernel $b$. This approach reduces the
considered problem to questions in harmonic analysis and lead us
with a complete answer.

\section{Preliminaries}
\setcounter{theorem}{0}
 \setcounter{equation}{0}

 Let $(\Omega,\mathcal{F},(\mathcal{F}_t)_{t\geq 0},P)$ be a complete
filtered probability space. By $T^d$ we denote the $d$-dimensional
torus which can be identified with the product $(-\pi,\pi]^d$. Let
$D(T^d)$ and $D'(T^d)$ denote, respectively, the 
topological vector space of functions 
$\phi\in C_0^\infty (T^d)$ endowed with the usual topology 
corresponding to the convergence (see, \cite{Ad} or \cite{Sc})
and the space of distributions. By
$\langle\xi,\phi\rangle$ we denote the value of a distribution
$\xi$ on a test function. We assume that $W_\Gamma$ is a
$D'(T^d)$-valued spatially homogeneous Wiener process (abbr.
s.h.~Wiener process) with covariance $\Gamma$. 
For more information about s.h.~Wiener process we refer to \cite{KaZa3}.

Any arbitrary s.h.~Wiener process $W_\Gamma$ is uniquely
determined by its covariance $\Gamma\in D'(t^d)$, which is 
a positive definite distribution, according to the formula
\begin{equation}{\label{eq2.1}}
 \mathbb{E} \langle W_\Gamma(t,\theta),\phi\rangle \,
 \langle W_\Gamma(\tau,\theta),\psi\rangle = \mathrm{min}(t,\tau)\,
 \langle \Gamma, \phi\star\psi_{(s)} \rangle \,,
\end{equation}
where $\phi,\psi\in D(T^d)$ and $\psi_{(s)}(\eta)=\psi(-\eta)$,
for $\eta\in T^d$. Because $W_\Gamma$ is spatially homogeneous
process, the distribution $\Gamma=\Gamma(\theta-\eta)$ for
$\theta,\eta\in T^d$.

The space covariance $\Gamma$ , like distribution in $D'(T^d)$,
may be uniquely expanded (see, e.g.\ \cite{GaWi} or \cite{Sc} )
into its Fourier series (with parameter $w=1$ because the period
is $2\pi$)
\begin{equation}{\label{eq2.2}}
\Gamma(\theta) = \sum_{n\in \mathbb{Z}^d} e^{i(n,\theta)}
\gamma_n, \quad \theta\in T^d,
\end{equation}
convergent in $D'(T^d)$. In (\ref{eq2.2}),
$(n,\theta)\!=\!\sum_{i=1}^d n_i\theta_i$ and $\mathbb{Z}^d$
denotes the product of integers.

The coefficients $\gamma_n$, in the Fourier series (\ref{eq2.2}),
satisfy (see, eg.\ \cite[Lesson 29, Section 29.4]{GaWi} 
or \cite{KaZa2}):
\begin{enumerate}
\item $\gamma_n=\gamma_{-n}$ for $n\in \mathbb{Z}^d$, $~~\gamma_n$ are 
non-negative,
\item the sequence $(\gamma_n)$ is {\em slowly increasing}, that is
\begin{equation}{\label{eq2.3}}
 \sum_{n\in \mathbb{Z}^d} \frac{\gamma_n}{1+|n|^r} < +\infty,
 \mbox{~~for~some~~} r>0.
\end{equation}
\end{enumerate}

Let us introduce, by induction, the following set of indexes.
Denote $\mathbb{Z}_s^1:=\mathbb{N}$, the set of natural numbers
and define
$\mathbb{Z}_s^{d+1}:=(\mathbb{Z}_s^1\times\mathbb{Z}^d)\cup
\{(0,n):n\in\mathbb{Z}_s^{d}\}$. Let us notice that
$\mathbb{Z}^d=\mathbb{Z}_s^d\cup(-\mathbb{Z}_s^d) \cup\{0\}$. For
instance, for $d=2$,
$\mathbb{Z}_s^2=\mathbb{N}\times\mathbb{Z}\cup\{(0,n):n\in\mathbb{Z}\}$.

Now, the s.h.\ Wiener process $W_\Gamma$ corresponding to the
covariance $\Gamma$ given by (\ref{eq2.2}) may be represented in
the form
\begin{eqnarray}{\label{eq2.4}}
 W_\Gamma(t,\theta) =
 \sqrt{\gamma_0}\beta_0(t) + \sum_{n\in \mathbb{Z}_s^d}
 \sqrt{2\gamma_n} \left[ \cos (n,\theta)\,\beta_n^1(t) \right.
  \!&\!+\!&\! \left.\sin (n,\theta)\,\beta_n^2(t) \right], \\ &&
 t\geq 0 \mbox{~and~} \theta\in T^d. \nonumber
\end{eqnarray}
In (\ref{eq2.4}) $\beta_0, \beta_n^1, \beta_n^2, n\in
\mathbb{Z}_s^d$, are independent real Brownian motions and
$\gamma_0, \gamma_n$ are coefficients of the series (\ref{eq2.2}).
The series (\ref{eq2.4}) is convergent in the sense of $D'(T^d)$.

Because any periodic distribution with positive period is a
tempered distribution (see, e.g.\ \cite{GaWi}), we may restrict our
considerations to the space $S'(T^d)$ of tempered distributions.
By $S(T^d)$ we denote the space of infinitely differentiable
rapidly decreasing functions on the torus $T^d$.

In the paper we use the following definition of the Sobolev spaces.
By $H^\alpha=H^\alpha(T^d)$, $\alpha\in\mathbb{R}$, we denote
the real Sobolev space of order $\alpha$ on the torus $T^d$. The
norms in such spaces may be expressed in terms of the Fourier
coefficients (see \cite{Ad})
\begin{equation}{\label{eq2.4a}}
||\xi||_{H^\alpha} \!= \!\left(\sum_{n\in\mathbb{Z}^d}(1+|n|^2)^\alpha
 |\xi_n|^2 \!\right)^\frac{1}{2} \!\!= \!\left(\! |\xi_0|^2 +
 2 \sum_{n\in \mathbb{Z}_s^d} (1+|n|^2)^\alpha \left( (\xi_n^1)^2
 +(\xi_n^2)^2 \right) \!\right)^\frac{1}{2} \!\!,
\end{equation}
where $\xi_n=\xi_n^1+i\xi_n^2,~\xi_=\bar{\xi}_{-n},~
n\in\mathbb{Z}^d$.

There is another possibility to define the Sobolev spaces (see,
e.g.\ \cite{ReSi}).

\section{Main results}
\setcounter{theorem}{0} \setcounter{equation}{0}

If $ b \in L^1_{loc} ( \mathbb{R}_+) $ and $ \mu \in \mathbb{C}$,
we shall denote by $ r(t,\mu) $ the unique solution in $
L^1_{loc}( \mathbb{R}_+ ) $ to the linear Volterra equation
\begin{equation}{\label{eq3.1}}
r(t,\mu) = b(t) + \mu \int_0^t b(t-s)r(s,\mu) ds, \quad t \geq 0,
\end{equation}
see \cite[Theorem 3.1]{GLS}. By taking Laplace transform in
$t$ we obtain $ r(\lambda, \mu) = \hat b(\lambda)(1-\mu \hat
b(\lambda))^{-1}$ and hence, in many cases, the function
$r(t,\mu)$ may be found explicitly. For example, if $ b(t)  \equiv
1, $ then $ r(\lambda ,\mu) = \frac{1}{\lambda -\mu}$ or,
inverting the Laplace transform, we obtain $ r(t,\mu) = e^{\mu t}.
$ Analogously, if $b(t) = e^{-t},$ then
 $ r(\lambda ,\mu) = \frac{1}{\lambda + (1- \mu)}$ and hence $ \quad r(t,\mu) = e^{
 (-1+\mu)t}$. Other examples are:
$$
\begin{array}{lclcl}
 & b(t) & = t, \quad & r(t,\mu)& = \displaystyle\frac{\sinh \sqrt{\mu}\,t}
{\sqrt{\mu}} \quad \mu \ne 0 \\  & b(t) & =
 te^{-t},  & \quad r(t,\mu)& =  e^{-t} \, \displaystyle \frac{\sinh
 \sqrt{\mu}\,t}{\sqrt{\mu}}, \quad \mu \neq 0.
\end{array}
$$
For more examples, see monograph \cite{Pr} by Pr\"uss.

\noindent Let us denote by $ \hat f(k), k \in \mathbb{Z}$, the
$k$-th Fourier coefficient of an integrable function~$f$:
$$ \hat f(k) = \frac{1}{2\pi} \int_0^{2\pi} e^{-ikt}f(t) dt. $$
Given $f$ defined on $[0,2\pi]$ we denote again by $f$ the
periodic extension to $\mathbb{R}$. Let $ b\in L^1 (
\mathbb{R}_+).$ We first observe from
$$F(t):= \int_{-\infty}^t
b(t-s) f(s) ds = \int_0^{\infty}b(s) f(t-s)ds $$ that $F$ is
periodic of period $T = 2\pi$ as $f$. Now using Fubini's theorem
we obtain,
\begin{equation}{\label{eq3.2}}
\hat F(k) = \tilde b(ik) \hat f(k), \quad k \in \mathbb{Z},
\end{equation}
where $ \tilde b(\lambda) := \displaystyle \int_0^{\infty}
e^{-\lambda t }\, b(t) dt $ denotes the Laplace transform of $b$.

In what follows we will assume that $\tilde b(ik)$ exists for all
$ k\in \mathbb{Z}$ and suppose that $ \lambda \to \tilde
b(\lambda) $ admits an analytical extension to a sector containing
the imaginary axis, and still denote this extension by $\tilde b$.
We introduce the following definition.

\begin{definition} \label{def3.2}
Let $r$ be given as in (\ref{eq3.1}). We say that a kernel $b\in
L^1 ( \mathbb{R_+})$ is admissible for equation (\ref{eq1.1}) if
$$ \lim_{|n| \to \infty} |n|^2 \int_{0}^{\infty}
[r(s, -|n|^2)]^2 ds =: C_b $$ exists.
\end{definition}

\begin{example}\label{ex3.3}

i) In case $ b(t) = e^{-t}$ we obtain $$|n|^2 \int_{0}^{\infty}
[r(s, -|n|^2)]^2 ds = \frac{|n|^2}{1+|n|^2} $$
and hence $C_b = 1.$\\

ii) In case $ b(t)= te^{-t}$ we obtain, after a calculation using
e.g. \cite[Formula 2.662 (2)]{GR},
$$|n|^2 \int_{0}^{\infty} [r(s, -|n|^2)]^2 ds = \frac{|n|^2}{4 +
4|n|^2}$$ and hence $ C_b = \frac{1}{4}.$
\end{example}

By a solution $X(t,\theta)$ to the stochastic Volterra equation
(\ref{eq1.1}) we will understand a process $X$ taking values in
the space $S'(T^d)$ and satisfying the equation (\ref{eq1.1}).
The following is our main result.

\begin{theorem}{\label{th3.4}}
Assume $b\in L^1( \mathbb{R_+}) $ is admissible for (\ref{eq1.1}).
Then, the equation (\ref{eq1.1}) has an $ H^{\alpha + 1}(
T^d)$-valued solution if and only if the Fourier coefficients
$(\gamma_n)$ of the covariance $\Gamma$ satisfy
\begin{equation}{\label{eq3.4}}
\sum_{n\in \mathbb{Z}^d} \gamma_n (1 + |n|^2)^{\alpha} < \infty.
\end{equation}
In particular, equation ~(\ref{eq1.1}) has an $L^2(T^d)$-valued
solution if and only if
$$ \sum_{n\in \mathbb{Z}^d} \frac{\gamma_n}{1+|n|^2}< \infty.$$
\end{theorem}

\bgproof We shall use the representation (\ref{eq2.4}) for the
Wiener process $W_\Gamma(t,\theta)$ with respect to the basis: $1,
\cos(n,\theta), \sin(n,\theta)$, where $n\in\mathbb{Z}_s^d$ and
$\theta\in T^d$. 
Denote by $W_n(t,\theta):=\cos (n,\theta)\,\beta_n^1(t)+\sin
(n,\theta)\,\beta_n^2(t),~ n\in\mathbb{Z}_s^d $, that is the
$n$-th element in the expansion (\ref{eq2.4}).
Equation (\ref{eq1.1}) may be solved
coordinatewise as follows.

Assume that
\begin{equation}\label{eq3.5}
 X(t,\theta)=\sum_{n\in\mathbb{Z}_s^d} [\cos(n,\theta)\,X_n^1(t) +
\sin(n,\theta)\,X_n^2(t)] + X_0(t).
\end{equation}

Introducing (\ref{eq3.5}) into (\ref{eq1.1}), we obtain
\begin{eqnarray*}
\cos(n,\theta)X_n^1(t) \!+\!\sin(n,\theta)X_n^2(t) \!\! &
\!\!=\!\!&\!\! -|n|^2\!\! \int_{-\infty}^t \!\!\!\!
b(t\!-\!s)[\cos(n,\theta)X_n^1(s) \! +\! \sin(n,\theta)X_n^2(s)] ds \\
\!\! &\!\!&\!\! + \!\!\sqrt{2\gamma_n} \int_{-\infty}^t\!\!\!\!
b(t\!-\!s) [\cos(n,\theta)\,\beta_n^1(s) \! +\!
\sin(n,\theta)\,\beta_n^2(s)] ds,
\end{eqnarray*}
or, equivalently
\begin{eqnarray*}
[\cos(n,\theta),\sin(n,\theta)] \left[\!\! \begin{array}{l}
X_n^1(t) \\ X_n^2(t) \end{array} \!\!\right] \!\! & \!\!=\!\!&\!\!
-|n|^2 \int_{-\infty}^t \!\!\!\! b(t-s)
[\cos(n,\theta),\sin(n,\theta)]
\left[\!\! \begin{array}{l} X_n^1(s) \\ X_n^2(s)\end{array}\!\!\right] ds \\
\!\! &\!+\!&\!\! \sqrt{2\gamma_n} \int_{-\infty}^t \!\!\!\!
b(t-s)[\cos(n,\theta),\sin(n,\theta)]
\left[\!\! \begin{array}{l} d\beta_n^1(s) \\
d\beta_n^2(s)\end{array}\!\!\right] .
\end{eqnarray*}
Denoting
$$ X_n(t):=
\left[ \begin{array}{l} X_n^1(t) \\ X_n^2(t)\end{array}\right]
\quad \mbox{and} \quad
 \beta_n(t):=
 \left[ \begin{array}{l} \beta_n^1(t) \\ \beta_n^2(t)\end{array}\right]
$$
we arrive at the equation
\begin{equation}{\label{eq3.6}}
 X_n(t) =-|n|^2 \int_{-\infty}^t b(t-s)X_n(s)ds + \sqrt{2\gamma_n}
\int_{-\infty}^t b(t-s) d \beta_n(s) \,.
\end{equation}
Taking Laplace transform in $t$, and making use of (\ref{eq3.1})
with $\mu = -|n|^2$ and (\ref{eq3.2}),
  we get the following solution to the equation (\ref{eq3.6}):
$$ X_n(t) = \int_{-\infty}^t r(t-s, -|n|^2) \sqrt{2 \gamma_n}
d \beta_n(s)= \int_0^{\infty} r(s, -|n|^2) \sqrt{2 \gamma_n} d
\beta_n(t-s).
$$
Hence, we deduce the following explicit formula for the solution
to the equation (\ref{eq1.1}):
\begin{equation}{\label{eq3.7}}
\begin{array}{rcl}
X(t,\theta) = \sqrt{\gamma_0} \beta_0(t) &+& \displaystyle \sum_{
n\in \mathbb{Z}_s^d } \sqrt{2 \gamma_n} \left[ \cos( n, \theta)
 \int_{0}^{\infty} r(s, -|n|^2) d\beta_n^1 (t-s) \right. \\ &+& \left.
\sin (n, \theta) \displaystyle \int_{0}^{\infty} r(s, -|n|^2) d
\beta_n^2 (t-s)\right].
\end{array}
\end{equation}
From definition of the space
$H^{\alpha}$ and condition (\ref{eq2.4a}) the process
 $ X(t) \in H^{\alpha + 1}, P$-almost surely, if and only if
\begin{equation}{\label{eq3.8}}
\sum_{n\in \mathbb{Z}^d} (1 + |n|^2)^{\alpha +1} \gamma_n \left[\!
\left(\!\! \int_{-\infty}^{t}\!\! r(t-s, -|n|^2) d \beta_n^1 (s)
\!\!\right)^2 \! + \! \left(\!\! \int_{-\infty}^{t}\!\! r(t-s,
-|n|^2) d \beta_n^2 (s) \!\!\right)^2 \right]\! < \!\infty.
\end{equation}
In the representation (\ref{eq2.4}) of $W_\Gamma$, processes 
$\beta_0, \beta_n^1,\beta_n^2, ~~n\in\mathbb{Z}_s^d,$ are
independent real Brownian motions.
So, the stochastic integrals in (\ref{eq3.8}) are independent
Gaussian random variables, too. Hence, (\ref{eq3.8}) holds 
$P$-almost surely if and only if
\begin{equation}{\label{eq3.9}}
\sum_{n\in \mathbb{Z}^d} (1 + |n|^2)^{\alpha +1} \gamma_n
\mathbb{E} [ (\int_{-\infty}^{t} r(t-s, -|n|^2) d \beta_n^1 (s))^2
+ (\int_{-\infty}^{t} r(t-s, -|n|^2) d \beta_n^2 (s))^2]< \infty.
\end{equation}
Or equivalently, using $L^2$ isometry of stochastic integrals, 
if and only if
\begin{equation}{\label{eq3.10}}
\sum_{n\in \mathbb{Z}^d} (1 + |n|^2)^{\alpha +1} \gamma_n
\int_{0}^{\infty} [ r(s, -|n|^2)]^2 ds < \infty.
\end{equation}
Since $b$ is admissible for the equation (\ref{eq1.1}), we
conclude that (\ref{eq3.10}) holds if and only if
\begin{equation*}
\sum_{n\in \mathbb{Z}^d} (1 + |n|^2)^{\alpha +1}
\frac{\gamma_n}{|n|^2} < \infty,
\end{equation*}
and the proof is complete.

 \edproof

 Concerning uniqueness, we have the following result.

\begin{proposition}
Assume $b\in L^1( \mathbb{R_+}) $ is admissible for (\ref{eq1.1})
and the following conditions hold:

(i)~ \quad $ \displaystyle \sum_{n\in \mathbb{Z}^d}
  \gamma_n (1 + |n|^2)^{\alpha} <\infty,$

 (ii) \quad $ \displaystyle \{1/\tilde b(ik) \}_
 {k\in \mathbb{Z}}
 \subset \mathbb{C} \setminus \{ -|n|^2 : n\in \mathbb{Z}^d \}$.\\

 Then, (\ref{eq1.1}) has a unique $ H^{\alpha + 1}(
T^d)$-valued solution.

\end{proposition}
\bgproof

Let $ X(t,\theta)$ be solution of $$ X(t,\theta) =
\int_{-\infty}^t b(t-s) \Delta X(s,\theta) ds.$$ Taking Fourier
transform in $\theta$ and denoting by $X_n(t)$ the $n$-th Fourier
coefficient of $X(t,\theta)$($t$ fixed), we obtain
$$ X_n(t) = -|n|^2 \int_{-\infty}^t b(t-s) X_n(s) ds $$ for all $
n \in \mathbb{Z}^d.$ Taking now Laplace transform in $t$, we get
that the Fourier coefficients of $X_n(t)$ ( $n$ fixed) satisfy
$$( 1 + |n|^2 \tilde b(ik) )\hat{X_n}(k)=0$$ for all $k \in
\mathbb{Z}$. According to (ii) we obtain that $ \hat{X_n}(k) = 0 $
for all $k \in \mathbb{Z}$ and all $ n \in \mathbb{Z}^d.$ Hence,
the assertion follows by uniqueness of the Fourier transform.

\edproof

 The following corollaries are an immediate consequence of Theorem \ref{th3.4}. The arguments are
  the same as in \cite{KaZa1}. We give here the
 proof for the sake of completeness.

 \begin{corollary}{\label{cor3.6}}
 Suppose $b\in L^1( \mathbb{R_+}) $ is admissible for (\ref{eq1.1}) and
 assume $\Gamma \in L^2 ( T^d). $ Then
 the integro-differential stochastic equation (\ref{eq1.1}) has a solution
 with values in $L^2( T^d)$ for $d=1,2,3.$

 \end{corollary}
 \bgproof
 We have to check equation (\ref{eq3.4}) with $\alpha = -1$. Note, that if
$\Gamma \in L^2 (T^d) $ then $\hat
 \Gamma = (\gamma_n) \in l^2 ( \mathbb{Z}^d).$ Consequently
\begin{equation*}
\sum_{n\in \mathbb{Z}^d} \frac{\gamma_n }{1 + |n|^2} \leq (
\sum_{n\in \mathbb{Z}^d} {\gamma_n }^2)^{1/2} ( \sum_{n\in
\mathbb{Z}^d} \frac{1}{(1 + |n|^2)^2})^{1/2}.
\end{equation*}
But $\sum_{n\in \mathbb{Z}^d} {\gamma_n }^2 < \infty $ and
\begin{equation}{\label{eq3.11}}
\sum_{n\in \mathbb{Z}^d} \frac{1}{(1 + |n|^2)^q} < \infty \mbox{
if and only if } 2q > d.
\end{equation}
Hence, the result follows.

 \edproof

 \begin{corollary}{\label{cor3.7}}
 Suppose $b\in L^1( \mathbb{R_+}) $ is admissible for (\ref{eq1.1}) and
 assume $\hat \Gamma \in l^p ( \mathbb{Z}^d) $ for $ 1 < p \leq 2.$ Then
 the integro-differential stochastic equation (\ref{eq1.1}) has a
 solution with values in $L^2( T^d)$ for all $ d <
 \frac{2p}{p-1}.$
\end{corollary}
\bgproof Note that
\begin{equation*}
\sum_{n\in \mathbb{Z}^d} \frac{\gamma_n }{1 + |n|^2} \leq (
\sum_{n\in \mathbb{Z}^d} {\gamma_n }^p)^{1/p} ( \sum_{n\in
\mathbb{Z}^d} \frac{1}{(1 + |n|^2)^q})^{1/q},
\end{equation*}
where $ \frac{1}{p} + \frac{1}{q} = 1.$ Hence the result follows
from (\ref{eq3.11}) with $q= \frac{p}{p-1}.$
\edproof

For $\alpha=-1$, the condition (\ref{eq3.4}) can be written as
follows.

\begin{theorem}{\label{th3.8}}
Let $b$ be admissible for the equation (\ref{eq1.1}). Assume that
the covariance $\Gamma$ is not only a positive definite
distribution but is also a non-negative measure. Then the equation
(\ref{eq1.1}) has $L^2(T^d)$-valued solution if and only if
\begin{equation}{\label{eq3.12}}
 (\Gamma,G_d) <+\infty,
\end{equation}
where
\begin{equation}{\label{eq3.13}}
 G_d(x) = \sum_{n\in\mathbb{Z}^d} \int_0^{+\infty}
 \frac{1}{\sqrt{(4\pi t)^d}}\, e^{-t}\, e^{-\frac{|x+2\pi n|^2}{t}}
 dt, \quad x\in T^d.
\end{equation}
\end{theorem}
The proof of Theorem \ref{th3.8} is the same that for Theorem 2,
part 2) in \cite{KaZa1}, so we omit it. 

The functions $G_d(x)$ are apparently well-known. If $d=1$, 
then $G_1$ is a positive continuous function. If $d=2$, 
$G_2$ is continuous outside 0 with singularity at 0 of the form
$\displaystyle c\,\log \frac{1}{|x|},~c>0$. 
For  $d>2$, $\displaystyle G_d(x)\sim c \frac{1}{|x|^{d-2}}$
as $x\rightarrow 0$. For more details
concerning the function $G_d$ we refer to \cite[Prop.~3]{KaZa2} and
\cite{La}.

From the above properties of functions $G_d$ defined by
(\ref{eq3.13}) and the condition (\ref{eq3.12}) we obtain the
following result (see \cite[Theorem 4]{KaZa2}).

\begin{corollary}{\label{cor3.9}}
Assume that $\Gamma$ is a non-negative measure and $b$ is
admissible. Then equation
(\ref{eq1.1}) has function valued solutions:\\
~~i) ~~\quad for all $\Gamma$ if $d=1$; \\
~ii) ~~\quad for exactly those $\Gamma$ for which
 $\int_{|\theta|\leq 1} \ln |\theta|\,\Gamma(d\theta)<+\infty$
 if $d=2$; \\
iii) \quad for exactly those $\Gamma$ for which
$\int_{|\theta|\leq 1} \frac{1}{|\theta|^{d-2}}\,
\Gamma(d\theta)<+\infty$ if $d\geq 3$.
\end{corollary}

In what follows, we will see that formula (\ref{eq3.7}) also
implies that $X$ is H\"older continuous with respect to $t$. In
order to do that, we need assumptions
very similar to those in \cite{ClDaPr}.\\

\noindent{\bf Hypothesis (H)}\\
{\it Assume that there exist $\delta\in (0,1)$ and $C_\delta>0$
such that, for all $s \in (-\infty,t)$ we have:\\
(i) $~\quad \displaystyle \int_s^t [r(t-\tau,-|n|^2)]^2 d\tau \leq
 C_\delta |n|^{2(\delta-1)}\,|t-s|^\delta$ ; \\
 (ii) $~\quad \displaystyle \int_{-\infty}^s [ r( t-\tau,-|n|^2)
 -r(s-\tau,-|n|^2)]^2 d\tau \leq C_\delta |n|^{2(\delta-1)}
 \,|t-s|^\delta$.}

\begin{proposition}
 Assume that $\displaystyle \sum_{n\in\mathbb{Z}^d}\frac{\gamma_n}
 {1+|n|^2}<+\infty$. Under Hypothesis (H), the trajectories of
 the solution $X$ to the equation
 (\ref{eq1.1}) are almost surely $\eta$-H\"older continuous
 with respect to $t$, for every $\eta\in (0,\delta/2)$.
\end{proposition}
\bgproof
 From the expansion (\ref{eq3.7}) and properties of stochastic integral,
 we have
 \begin{eqnarray*}
 \displaystyle
  \mathbb{E} \!\! & \!\! ||X(t,\theta)-X(s,\theta)||^2_{L^2} =
 \mathbb{E}\left|\left| \sqrt{\gamma_0} (\beta_0(t)-\beta_0(s))
  \hspace{38ex} \right.\right. \\
 & \!\!+\! \displaystyle\sum_{n\in\mathbb{Z}_s^d} \sqrt{2\gamma_n} \left[
 \cos (n,\theta)\!\left(
 \displaystyle \int_{-\infty}^t r(t\!-\!\tau,\!-\!|n|^2) d\beta_n^1(\tau) \!
 -\! \displaystyle\int_{-\infty}^s r(s\!-\!\tau,\!-\!|n|^2) d\beta_n^1(\tau)\right)
 \right.\\
 & \hspace{4ex} \left.\left.
 +\sin (n,\theta)\left( \left.
  \displaystyle\int_{-\infty}^t r(t\!-\!\tau,\!-\!|n|^2) d\beta_n^2(\tau)
 - \displaystyle\int_{-\infty}^s r(s\!-\!\tau,\!-\!|n|^2) d\beta_n^2(\tau)
 \right) \right] \right|\right|^2_{L^2} \\
 = \!\! & \!\! (2\pi)^d \left(\gamma_0 |t-s| \right. \hspace{66ex} \\
 \!\! & \!\! \left. + \displaystyle \sum_{n\in\mathbb{Z}_s^d}
 2\gamma_n\left[  \displaystyle\int_{-\infty}^s [r(t\!-\!\tau,\!-\!|n|^2)
 - r(s\!-\!\tau,\!-\!|n|^2)]^2d\tau +  \displaystyle\int_s^t r^2(t\!-\!\tau,\!-\!|n|^2)d\tau
 \right] \right).
 \end{eqnarray*}

According to assumptions (i) and (ii) of the Hypothesis (H), we
have
$$ \mathbb{E}||X(t,\theta)-X(s,\theta)||^2_{L^2} \leq C_\delta
\sum_{n\in\mathbb{Z}_s^d} 2\gamma_n \,|n|^{2(\delta-1)}
 |t-s|^\delta . $$
Because $X$ is a Gaussian process, then for any $m\in N$, there
exists a constant $C_m>0$ that
$$\mathbb{E}||X(t,\theta)-X(s,\theta)||^{2m}_{L^2} \leq C_m
\left[ C_\delta \sum_{n\in\mathbb{Z}_s^d} 2\gamma_n
\,|n|^{2(\delta-1)} \right]^m |t-s|^{m\delta} . $$ 
Taking $m$ such
that $m\delta>1$ and using the Kolmogorov's criterion for
continuity (see, e.g.\ \cite{Str}), we obtain that the solution
$X(t,\theta)$ is $\eta$-H\"older continuous, with respect to $t$,
for $\eta=\delta/2-1/(2m)$. \edproof
\begin{example}
 Let us consider the particular case $b(t) = e^{-t}, ~t\geq 0$.
 Then, by previous considerations,
 $r(t,-|n|^2) = e^{(-1-|n|^2)\,t}$. One can check that in this
 case the Hypothesis (H) is fulfilled.
\end{example}

\begin{remark} We observe that the condition (i) in Hypothesis (H)
is the same as
$$
|n|^2 \int_0^{t} [r(s, -|n|^2)]^2 ds \leq C_{\delta} |n|^{2\delta}
|t|^{\delta}
$$
and hence it is nearly equivalent to say that the function $b$ is
admissible.
\end{remark}

{\bf Acknowledgement} The authors would like to thank the
anonymous referee for the careful reading of the manuscript. 
The valuable remarks made numerous improvements throughout.

\end{document}